\documentclass[a4paper,11pt]{article}
\usepackage[utf8]{inputenc}
\usepackage[T1]{fontenc}
\usepackage{lmodern}
\usepackage{microtype}
\usepackage[a4paper,textwidth=16cm,textheight=22cm,centering]{geometry}

\usepackage{mathtools}
\usepackage{amssymb,amsfonts}
\usepackage{amsthm}
\usepackage{setspace}
\usepackage{indentfirst}

\usepackage{cite}
\usepackage{xcolor}
\usepackage[ ]{hyperref}
\hypersetup{
	colorlinks=true,   
	linktoc=all,        
	linkcolor=red,     
	citecolor=blue,
	urlcolor=black,} 

\theoremstyle{plain}
\newtheorem{theorem}{Theorem}[section]
\newtheorem{lemma}[theorem]{Lemma}
\newtheorem{claim}[theorem]{Claim}
\newtheorem{fact}[theorem]{Fact}
\newtheorem{corollary}[theorem]{Corollary}
\newtheorem{proposition}[theorem]{Proposition}

\theoremstyle{definition}
\newtheorem{definition}[theorem]{Definition}
\newtheorem{problem}[theorem]{Problem}

\theoremstyle{remark}

\newcommand{\bhz}[1]{\textcolor{black}{#1}}

\DeclareMathOperator{\twr}{twr}
\begin{document}
	\date{}
	\begin{spacing}{1.03}
\title{Hypergraph Ramsey Numbers via Oriented Palette Constructions}
\author{Hanzhi Bai$^1$, \;\;\; Xinyu Hu$^2$, \;\;\; Kaidi Liu$^1$, \;\;\; Guanghui Wang$^{1}$}

\footnotetext[1]{{School of Mathematics, Shandong University, Jinan 250100, P.~R.~China.
Emails: {\tt bhz@mail.sdu.edu.cn} (H. Bai), {\tt 202611923@mail.sdu.edu.cn} (K. Liu), {\tt ghwang@sdu.edu.cn} (G. Wang). Supported by the National Natural Science Foundation of China (12231018) and the State Key Laboratory of Cryptography and Digital Economy Security.}}
\footnotetext[2]{Data Science Institute, Shandong University, Jinan 250100, P.~R.~China.
Email: {\tt huxinyu@sdu.edu.cn} (X. Hu).
Supported by National Postdoctoral Fellowship Program (C-tier) (GZC20252005) and NSFC (No. 12601677).}

\maketitle	

\begin{abstract}
Given a family $\mathcal F$ of $k$-graphs, the Ramsey number $r_k(\mathcal F,n)$ is the minimum $N$ such that every red/blue coloring of the edges of $K_N^{(k)}$ yields either a red copy of some $F\in\mathcal F$ or a blue copy of $K_n^{(k)}$. In particular, if $\mathcal F=\{F\}$, we write $r_k(F,n)$. 

In this paper, we use oriented palette constructions and the hypergraph container method to obtain a lower bound on $r_k(\mathcal F,n)$ for special families $\mathcal F$. More precisely, we independently label lower-dimensional faces and use the palette to exclude prescribed red configurations, while a container argument controls the size of blue cliques.

As a first application of our theorem, we obtain \(r_k(B^{(k)},n)>2^{\Omega_k(n\log n)}\). Combined with the upper bound due to Mubayi and Suk, this yields \(r_k(B^{(k)},n)=2^{\Theta_k(n\log n)}\), where \(B^{(k)}\) is the half-graph. This answers in the affirmative a question posed by Mubayi and Suk. As a second application of our theorem, we show $r_k(k+1,3;n)=2^{\Theta_k(n^{k-2}\log n)}$, which, together with a new stepping-up lemma developed in this paper, further yields $r_k(k+1,t;n)=\twr_{t-1}(\Theta_{k,t}(n^{k-t+1}\log n))$ for every fixed $k\ge6$ and $4\le t\le k-2$, where $\twr_1(x)=x$ and $\twr_{j+1}(x)=2^{\twr_j(x)}$. Here $r_k(k+1,t;n)$ denotes the minimum $N$ such that every red/blue coloring of the $k$-subsets of an $N$-vertex set contains either $k+1$ vertices spanning at least $t$ red edges or a blue clique of order $n$. As a third application of our theorem, we show that $r_k(F_{\ell,k,d},n)>2^{\Omega_{k,\ell,d}(n^{k-d-2}\log n)}$ for all fixed integers $k\ge4$, $1\le d\le k-3$, and $\ell\ge2k+1$ with $k\nmid\ell$, where $F_{\ell,k,d}$ is obtained from the $k$-uniform tight cycle of length $\ell$ by deleting $d$ consecutive edges.
\end{abstract}

\section{Introduction}\label{sec:introduction}
A $k$-uniform hypergraph, or $k$-graph, consists of a vertex set and a collection of $k$-element subsets called edges. Write $K_n^{(k)}$ for the complete $k$-graph on $n$ vertices. Given a family $\mathcal F$ of $k$-graphs, the Ramsey number $r_k(\mathcal F,n)$ is the minimum $N$ such that every red/blue coloring of the edges of $K_N^{(k)}$ yields either a red copy of some $F\in\mathcal F$ or a blue copy of $K_n^{(k)}$. In particular, if $\mathcal F=\{F\}$, we write $r_k(F,n)$. In this paper, we focus on the hypergraph Ramsey problem ($k\ge 3$). For the graph case, we refer the reader to Morris \cite{R.morris} for a nice survey on this topic.

The quantitative study of diagonal hypergraph Ramsey numbers goes back to Erd\H{o}s and Rado~\cite{ER} and Erd\H{o}s, Hajnal, and Rado~\cite{EHR}. Specifically, they showed that $\operatorname{twr}_{k-1}(\Omega(n^2))< r_k(n,n)< \operatorname{twr}_{k}(O(n))$, where the tower function is defined by $\operatorname{twr}_{1}(x)=x$ and $\operatorname{twr}_{i+1}(x)=2^{\operatorname{twr}_{i}(x)}$ for $i\ge 1$. It is a major open problem to determine if $r_k(n,n)> \operatorname{twr}_{k}(\Omega(n))$. For the off-diagonal hypergraph Ramsey numbers $r_k(s,n)$, Erd\H{o}s and Hajnal \cite{EH} conjectured that for every fixed $s>k\ge 4$, one has $r_k(s,n)>\operatorname{twr}_{k-1}(\Omega(n))$. This conjecture was independently verified for $k\ge4$ and $s\ge k+3$ by Mubayi and Suk \cite{MSoffdiag} and by Conlon, Fox and Sudakov \bhz{\cite{CFS}}. Recently, Du, Hu, Liu and Wang \cite{D-H-L-W-4} verified it for $k\ge 4$ and $s=k+2$. For the remaining case $s=k+1$, Du, Hu, Liu and Wang \cite{D-H-L-W-3} obtained the correct tower height.

For a non-complete $k$-graph $F$, there are also many results on $r_k(F,n)$~\cite{CFGHMSV,CFGHMSVY,MR}. Let $K_4^{(3)}\setminus e$ denote the $3$-graph on four vertices obtained by removing one edge from $K_4^{(3)}$. Erd\H{o}s and Hajnal~\cite{EH} proved that $r_3(K_4^{(3)}\setminus e,n)<(n!)^2$. Fox and He~\cite[Theorem~1.1]{FH} showed that $r_3(K_4^{(3)}\setminus e,n)=2^{\Theta(n\log n)}$. Mubayi and Suk \cite{MSconstruction} generalized $K_4^{(3)}\setminus e$ to the \emph{$k$-half-graph} $B^{(k)}$ for $k\ge 3$.

For $k\ge3$, choose disjoint \bhz{vertex} sets $S,T$ of size $k-1$ and a vertex $u\in S$. The \emph{$k$-half-graph} $B^{(k)}$ has edge set
\begin{equation}\label{eq:half-graph}
 E(B^{(k)})=\{S\cup\{t\}:t\in T\}\cup\{T\cup\{u\}\}.
\end{equation}
Thus $B^{(3)}=K_4^{(3)}\setminus e$. Mubayi and Suk~\cite[Theorems~4.1 and~4.3]{MSconstruction} showed that 
\begin{align}\label{half-up-lw}
    2^{c_kn}<r_k(B^{(k)},n)\le(n!)^{k-1}
\end{align}
for every fixed $k\ge4$ and all sufficiently large $n$. These estimates leave a factor of $\log n$ between the lower and upper exponents. Mubayi and Suk posed the following problem in their survey~\cite[Problem~8.1]{MSsurvey}.

\bhz{\begin{problem}[Mubayi and Suk \cite{MSsurvey}]\label{prob:MS}
For each $k\ge3$, does there exist $c_k>0$ such that $r_k(B^{(k)},n)>2^{c_kn\log n}$ for all sufficiently large $n$?
\end{problem}}

The generalized Ramsey number \(r_k(s,t;n)\), introduced by Erd\H{o}s and Hajnal \cite{EH}, is defined as the smallest positive integer \(N\) such that every red/blue coloring of the edges of \(K_N^{(k)}\) contains either an \(s\)-vertex set spanning at least \(t\) red edges or a blue \(K_n^{(k)}\). \bhz{In particular,}
\(
r_k\left(s,\binom{s}{k};n\right)=r_k(s,n).
\)

For fixed \(k\ge3\), R\"{o}dl and \v{S}inajov\'a \cite{Rodl-sinajova} (for partial Steiner systems), along with Kostochka, Mubayi, and Verstra\"ete \cite{K-M-V} (for independent sets in hypergraphs), showed that
\[
\Omega(n^{k-1}/\log n) < r_k(k+1, 2; n) < O( n^{k-1}/\log n).
\] Mubayi and Suk \cite[Theorem~2]{MSEH} proved that
\begin{align}\label{r_k(k+1,3;n)}
     2^{c_kn^{k-2}}< r_k(k+1, 3; n)< 2^{C_kn^{k-2}\log n}.
\end{align}
For higher local parameters, they \cite[Theorem~3]{MSEH} used the vertex-online ordered Ramsey game and the stepping-up lemma to obtain, for \(4 \le t \le k - 2\),
\begin{align}\label{r_k(k+1,t;n)}
     \operatorname{twr}_{t-1}(O(n^{k-t+1}\log n)) > r_k(k+1, t; n) > 
\begin{cases} 
\operatorname{twr}_{t-1}(\Omega(n^{k-t+1})) & \text{if } k - t \text{ is even,} \\
\operatorname{twr}_{t-1}(\Omega(n^{(k-t+1)/2})) & \text{if } k - t \text{ is odd.}
\end{cases}
\end{align}
These results settle the  Erd\H{o}s-Hajnal conjecture \cite{EH} for $r_k(s, t; n)$ in almost all cases with $s=k+1$. Following this, Mubayi, Suk and Zhu \cite{M-S-Z} showed that
\[
r_k(k+1,k-1;n)=\operatorname{twr}_{k-2}(n^{\Theta(1)}).
\]
For \(k\ge4\), in the remaining case \(r_k(k+1,k;n)\), the best known result is
\[
r_k(k+1,k;n)> \operatorname{twr}_{k-2}(n^{\Omega(\log \log n)}),
\]
which was recently obtained by Fan, Hu, Lin and Lu \cite{F-H-L-L}.

For $s\ge k+2$, the growth of $r_k(s,t;n)$ exhibits a polynomial-to-exponential transition as $t$ increases: it is polynomial in $n$ for sufficiently small $t$, while for sufficiently large $t$ it admits an exponential lower bound in a positive power of $n$. We refer the reader to~\cite{EH,MR,A-H-Y} for details.
\medskip

In this paper, we answer Problem \ref{prob:MS} affirmatively by proving the following theorem.
\begin{theorem}\label{thm:half}
For every fixed \(k\ge3\), we have
\[
r_k(B^{(k)},n)>2^{\Omega_k(n\log n)},
\]
which, together with the upper bound from (\ref{half-up-lw}), yields
\[
r_k(B^{(k)},n)=2^{\Theta_k(n\log n)}.
\]
\end{theorem}

On the other hand, we improve the lower bound in (\ref{r_k(k+1,3;n)}) and obtain
\begin{theorem}\label{thm:localbase}
For every fixed \(k\ge3\), we have
\[
r_k(k+1,3;n)=2^{\Theta_k(n^{k-2}\log n)}.
\]
\end{theorem}
Combining Theorem \ref{thm:localbase} with the proof in \cite[Theorem~3]{MSEH}, we obtain that, for \(4 \le t \le k - 2\) with \(k-t\) even,
\begin{align}\label{mubayi-suk results}
    r_k(k+1, t; n)=\operatorname{twr}_{t-1}(\Theta(n^{k-t+1}\log n)).
\end{align}

For \(4 \le t \le k - 2\) with \(k-t\) odd, we obtain the same result by proving the following new stepping-up lemma.

\begin{lemma}\label{su:oddstep}
For every odd integer \(k\ge7\) and integer \(n\ge k\),
\begin{equation}\label{eq:oddstep-main}
 r_k(k+1,4;2kn)>2^{r_{k-1}(k,3;n)-1}.
\end{equation}
\end{lemma}

Consequently, we prove the following theorem, in which the upper bound is obtained from (\ref{r_k(k+1,t;n)}).
\begin{theorem}\label{su:tower}
For fixed integers \(k\ge6\) and \(4\le t\le k-2\), there exist positive constants \(c_{k,t},C_{k,t}\) such that, for all sufficiently large \(n\),
\[
\operatorname{twr}_{t-1}(c_{k,t}n^{k-t+1}\log n)
 < r_k(k+1, t; n)
 <\operatorname{twr}_{t-1}(C_{k,t}n^{k-t+1}\log n).
\]
\end{theorem}

Indeed, we establish a more general result in this paper. We state the relevant definitions and theorems below. 

\begin{definition}\label{Delta_k,a}
For integers \(k\) and \(a\) with \(1\le a\le\lfloor k/2\rfloor\), define the \(k\)-graph \(\Delta_{k,a}\) on the vertex set \(V = U \cup X \cup Y \cup Z\), with edge set
\[
E(\Delta_{k,a}) = \{U\cup X\cup Y,\ U\cup X\cup Z,\ U\cup Y\cup Z\},
\]
where \(U, X, Y, Z\) are pairwise disjoint, \(|U|=k-2a\), and \(|X|=|Y|=|Z|=a\). Consequently, \(\Delta_{k,a}\) has pairwise edge intersections of size \(k-a\).
\end{definition}

\begin{definition}\label{T_k,r}
For integers \(k\) and \(r\) with \(2\le r\le k-2\), define the \(k\)-graph \(T_{k,r}\) on the vertex set \(V = S \cup \{x,y\} \cup W\), with edge set
\begin{equation}\label{eq:Tdefinition}
E(T_{k,r}) = \left\{ S\cup\{x\},\ S\cup\{y\},\ U\cup\{x,y\}\cup W \right\},
\end{equation}
where \(S, \{x,y\}, W\) are pairwise disjoint, \(|S|=k-1\), \(|W|=k-r-1\), and \(U\subseteq S\) is some fixed subset of size \(r-1\).
\end{definition}

For $2\le r<k$, define a family of $k$-graphs
\begin{equation}\label{eq:family}
 \mathcal F_{k,r}=
 \begin{cases}
 \{\Delta_{k,a}:1\le a\le\min(k-r,\lfloor k/2\rfloor)\}\cup\{T_{k,r}\},&r\le k-2,\\
 \{\Delta_{k,1}\},&r=k-1.
 \end{cases}
\end{equation}
\bhz{Each} member has three edges, but the family generally contains nonisomorphic hypergraphs.   

The same definition of $T_{k,r}$ can be extended to $r=k-1$ by taking $W=\varnothing$. In that case, writing $S\setminus U=\{z\}$ gives the three edges $U\cup\{x,z\}$, $U\cup\{y,z\}$, and $U\cup\{x,y\}$. Thus $T_{k,k-1}\cong\Delta_{k,1}$, which explains its omission as a separate member of $\mathcal F_{k,k-1}$.

\begin{theorem}\label{thm:family}
Fix integers $k\ge3$ and $2\le r<k$. There is a constant $c_{k,r}>0$ such that, for all sufficiently large $n$,
\[
 r_k(\mathcal F_{k,r},n)> 2^{c_{k,r}n^{r-1}\log n}.
\]
\end{theorem}
The conclusion uses a single coloring avoiding every member of $\mathcal F_{k,r}$. It does not follow by combining separate lower bounds for the members. Choosing $r=k-a$ gives the following consequence.
\begin{corollary}\label{cor:delta}
For every fixed $k\ge3$ and $1\le a\le\lfloor k/2\rfloor$, there is $c_{k,a}>0$ such that $r_k(\Delta_{k,a},n)> 2^{c_{k,a}n^{k-a-1}\log n}$ for all sufficiently large $n$.
\end{corollary}

For $k\ge4$, the half-graph $B^{(k)}$ contains $T_{k,2}$. For $k=3$, $B^{(3)}=\Delta_{3,1}$. Consequently, Theorem~\ref{thm:half} is a direct consequence of Theorem~\ref{thm:family}. On the other hand, any three distinct $k$-edges on a $(k+1)$-set form a copy of $\Delta_{k,1}$, so $r_k(\Delta_{k,1},n)=r_k(k+1,3;n)$. The choice $r=k-1$ in Theorem~\ref{thm:family} therefore yields Theorem \ref{thm:localbase}.
\medskip

Conlon, Fox, Gunby, He, Mubayi, Suk, Verstra\"ete, and Yu~\cite[Conjecture~1.1]{CFGHMSVY} conjectured that a fixed $3$-graph $H$ satisfies $r_3(H,n)\le n^{c(H)}$ for some $c(H)>0$ and all sufficiently large $n$ if and only if $H$ is \emph{iterated tripartite}, that is, contained in an iterated blowup of a $3$-edge obtained by recursively inserting complete tripartite $3$-graphs into its parts.

Recently, Conlon, He, Nie, and Post~\cite{CHNP} verified this conjecture for $3$-graphs with at most five tight components and reduced the general case to linear $3$-graphs. Subsequently, Gishboliner, Hajebi, and Li~\cite{GHL} settled the conjecture in full.

For $k\ge4$, the analogous classification remains open. In an iterated blowup of a $k$-edge, any subhypergraph whose edges are connected through overlaps of at least two vertices is $k$-partite. Motivated by this observation, we use the same palette construction to exclude entire families of non-$k$-partite $k$-graphs with strong overlap connectivity, even with arbitrarily many edges, and obtain explicit exponential Ramsey lower bounds.

A $k$-graph is \emph{$k$-partite} if its vertices admit a partition into $k$ classes such that every edge meets each class exactly once. A vertex is \emph{isolated} if it belongs to no edge. For a $k$-graph $F$ and $s\in[k]$, let $\Gamma_s(F)$ be the graph on $E(F)$ in which distinct edges are adjacent if they share at least $s$ vertices, and let $\Gamma_s(F;v)$ be its induced subgraph on the edges containing $v$. A one-vertex graph is considered connected. For $k\ge4$ and $2\le r\le k-2$, let $\mathcal H_{k,r}$ be the family of non-$k$-partite $k$-graphs $F$ without isolated vertices satisfying
\begin{enumerate}
\renewcommand{\labelenumi}{(\roman{enumi})}
\item $\Gamma_{r+1}(F)$ is connected;
\item $\Gamma_r(F;v)$ is connected for every $v\in V(F)$.
\end{enumerate}
\medskip

{\color{black}\noindent \textit{Remark.}~
The two connectivity conditions play different roles: connectedness of $\Gamma_{r+1}(F)$ forces all red edges to have the same palette support, while connectedness of each $\Gamma_r(F;v)$ fixes the position of $v$ across its incident red edges. Proposition~\ref{oc:criterion} therefore makes any red copy of $F$ $k$-partite, contradicting the definition of $\mathcal H_{k,r}$. This gives the following theorem.}

\begin{theorem}\label{oc:ramsey}
Fix integers $k\ge4$ and $2\le r\le k-2$. There exists $c_{k,r}>0$ such that, for all sufficiently large $n$,
\[
 r_k(\mathcal F_{k,r}\cup\mathcal H_{k,r},n)
 >2^{c_{k,r}n^{r-1}\log n}.
\]
\end{theorem}
In particular, this bound holds for $r_k(\mathcal H_{k,r},n)$ and for every fixed member of $\mathcal H_{k,r}$. The theorem asserts simultaneous avoidance by one coloring, with a constant independent of the choice of a member. The two connectivity conditions are sufficient, but they do not include all configurations in $\mathcal F_{k,r}$.

Conlon, Fox, Gunby, He, Mubayi, Suk, Verstra\"ete and Yu~\cite[Theorem~4.1]{CFGHMSVY} proved $r_k(F,n)>2^{\Omega_F(n^{2/k})}$ whenever $F$ is non-$k$-partite and $\Gamma_2(F)$ is connected. For $F\in\mathcal H_{k,r}$, Theorem~\ref{oc:ramsey} improves the exponent to $n^{r-1}\log n$, with one coloring avoiding the entire family and a constant depending only on $k,r$. It does not establish the higher-tower predictions in~\cite[Conjectures~4.2 and~4.3]{CFGHMSVY}.

To give an application, let $C_\ell^{(k)}$ be the tight cycle on vertices $v_0,\ldots,v_{\ell-1}$, whose edges are $E_i=\{v_i,v_{i+1},\ldots,v_{i+k-1}\}$ for $0\le i<\ell$, with subscripts interpreted modulo $\ell$. Let $F_{\ell,k,d}$ be obtained by deleting $E_0,\ldots,E_{d-1}$.
\begin{corollary}\label{oc:cycles}
Fix integers $k\ge4$, $1\le d\le k-3$, and $\ell\ge2k+1$ with $k\nmid\ell$. There exists $c_{k,\ell,d}>0$ such that, for all sufficiently large $n$,
\[
 r_k(F_{\ell,k,d},n)>2^{c_{k,\ell,d}n^{k-d-2}\log n}.
\]
In particular, $r_k(C_\ell^{(k)}\setminus e,n)>2^{\Omega_{k,\ell}(n^{k-3}\log n)}$.
\end{corollary}
We prove the corollary by showing that $F_{\ell,k,d}\in\mathcal H_{k,k-d-1}$. In the stated range, the condition $k\nmid\ell$ is exactly the obstruction to $k$-partiteness. {\color{black}For comparison, the ordinary (unordered) $k$-uniform tight path $P_s^{(k)}$ satisfies $r_k(P_{k+1}^{(k)},n)=\Theta_k(n^{k-1}/\log n)$ and 
\[\Omega_{k,s}(n^{k-1}/\log n)\le r_k(P_s^{(k)},n)\le O_{k,s}(n^{k-1})\]
for fixed $s\ge k+2$ \cite{Rodl-sinajova,K-M-V}. For $k=3$, results are known: $r_3(P_4^{(3)},n)=\Theta(n^2/\log n)$, $\Omega(n^2/\log n)\le r_3(P_5^{(3)},n)\le O(n^2)$, and $r_3(P_s^{(3)},n)=\Theta_s(n^2)$ for $s\ge 6$ \cite{Phelps-Rodl,Cooper-Mubayi,Mattheus-Mubayi-Nie-Verstraete}. Although $F_{\ell,k,d}$ contains $P_{\ell-d}^{(k)}$, Corollary~\ref{oc:cycles} yields an exponential Ramsey lower bound, in contrast to the polynomial growth for ordinary tight paths.} The resulting bound is single exponential and is likely not optimal. The range $k\ge4$ does not include the question about $C_\ell^{(3)}\setminus e$ discussed in~\cite[Section~4]{CFGHMSVY}.
\medskip

The palette construction develops the use of pair colors in Conlon, Fox, Gunby, He, Mubayi, Suk, and Verstra\"ete~\cite[Section~2.2]{CFGHMSV}. The auxiliary counting hypergraph is also motivated by Lin and Niu~\cite[Sections~3.3 and~3.4]{QL}. We use the container method developed independently by Balogh, Morris, and Samotij~\cite{BMS} and Saxton and Thomason~\cite{ST}. Our palette-based bounds are single exponential and do not establish the higher-tower predictions in~\cite[Conjectures~4.2 and~4.3]{CFGHMSVY}.
\medskip

\noindent\textit{Organization}. Section~\ref{sec:preliminaries} introduces the notation and preliminary results. Section~\ref{sec:construction} constructs the palettes, and Section~\ref{sec:avoidance} proves the deterministic exclusion results. Section~\ref{sec:counting} establishes the counting estimate. Finally, Section~\ref{sec:main-proof} derives the Ramsey lower bounds and proves the odd-uniformity stepping-up lemma and its tower consequences.

\section{Terminology and Preliminary Results}\label{sec:preliminaries}
\subsection{Notation}
For two positive integers $s<t$, let $[s]=\{1,\ldots,s\}$ and $[s,t]=\{s,s+1,\ldots,t\}$. Let $\binom Xs$ be the family of $s$-element subsets of a set $X$. A $k$-graph $H$ has vertex set $V(H)$ and edge set $E(H)\subseteq\binom{V(H)}k$. We write $v(H)=|V(H)|$ and $e(H)=|E(H)|$. For $X\subseteq V(H)$, the induced hypergraph $H[X]$ has edge set $E(H)\cap\binom Xk$. A copy of $F$ in $H$ is an injective map of vertices that takes each edge of $F$ to an edge of $H$. The hypergraph $H$ is $\mathcal F$-free if it has no copy of any member of $\mathcal F$. A set $X$ is independent if $e(H[X])=0$, and $\alpha(H)$ is its maximum possible size. Thus $r_k(\mathcal F,n)$ is the least $N$ such that every $\mathcal F$-free $k$-graph of order $N$ has independence number at least $n$. An independent set in the red hypergraph is also called a blue set.

\subsection{Palettes and face labels}
A hypergraph is \emph{linear} if two distinct edges intersect in at most one vertex. A \emph{Berge triangle} consists of distinct edges $A_1,A_2,A_3$ and distinct vertices $x_{12},x_{23},x_{31}$ such that $x_{ij}\in A_i\cap A_j$. No additional restriction on intersections is imposed.

Fix $2\le r<k$ and put $m=\binom kr$. The $r$-subsets of a $k$-set are called its $r$-faces. If $E=\{v_1<\cdots<v_k\}$, list its $r$-faces $\sigma(E)$ as $\sigma_1(E),\ldots,\sigma_m(E)$ in lexicographic order of their increasing lists of positions. 

A \emph{palette} is a set of vectors. An \emph{oriented palette} is a set $\mathcal P\subseteq[q]^m$ of vectors with distinct coordinates, with at most one vector for each coordinate support. The support of $\mathbf a=(a_1,\ldots,a_m)$ is $\{a_1,\ldots,a_m\}$. Write $\mathcal A(\mathcal P)$ for the $m$-graph on $[q]$ consisting of these supports. We write $\langle v_1,\ldots,v_k\rangle$ for a $k$-set $\{v_1,\ldots,v_k\}$ with $v_1<\cdots<v_k$. Given a face-labeling map $f:\binom{[N]}r\to[q]$, define $R(f,\mathcal P)$ to be the red 
$k$-graph on $[N]$ whose edges consist of the $k$-sets $E$ satisfying
\begin{equation}\label{eq:face-rule}
 \bigl(f(\sigma_1(E)),\ldots,f(\sigma_m(E))\bigr)\in\mathcal P.
\end{equation}
The other $k$-sets are blue. Palette coordinates are indexed by $r$-subsets of positions in a $k$-set, whereas vertices of $\mathcal A(\mathcal P)$ are colors. These are separate from the vertices of $R(f,\mathcal P)$.

\subsection{Supersaturation}
We also use the following standard supersaturation consequence of the Tur\'{a}n theorem for complete hypergraphs. Supersaturation results strengthen extremal results by showing that once the edge density exceeds the Tur\'{a}n density by a fixed positive amount, many copies of the forbidden hypergraph must occur. This phenomenon was studied systematically by Erd\H{o}s and Simonovits \cite{ES} for graphs and uniform hypergraphs. See also Keevash \cite[Lemma 2.1]{Keevash} for a convenient formulation in terms of Tur\'{a}n density. In particular, for every fixed $r$-graph $F$ and every $\delta >0$, there exists a constant $c=c(F,\delta)$ such that every sufficiently large $r$-graph $J$ on $n$ vertices with $e(J)\ge (\pi (F)+\delta)\binom{n}{r}$ contains at least $cn^{v(F)}$ copies of $F$. We apply this result to complete $r$-graphs. Recall that \(\pi(K_k^{(r)})=\lim_{n\to\infty}\operatorname{ex}(n,K_k^{(r)})/\binom{n}{r}\). De Caen\cite{decaen} proved \(\pi(K_k^{(r)})\le 1-1/\binom{k-1}{r-1}\). Combining this bound with the supersaturation theorem above gives the following consequence.
\begin{lemma}\label{lem:supersaturation}
For all integers \(k>r\ge 2\) and every \(\varepsilon>0\), there exists a constant \(\kappa=\kappa(k,r,\varepsilon)>0\) such that, for all sufficiently large \(n\), every \(r\)-graph \(J\) on \(n\) vertices with
\[e(J)\ge\left(1-1/\binom{k-1}{r-1}+\varepsilon\right)\binom{n}{r}\]
contains at least \(\kappa n^k\) copies of \(K_k^{(r)}\).
\end{lemma}

\subsection{Hypergraph containers}
For an $m$-graph $G$, its average degree is $d(G)=m e(G)/v(G)$. For $2\le j\le m$, let $\Delta_j(G)$ denote the maximum number of edges containing a fixed set of $j$ vertices. The hypergraph container method of Balogh, Morris and Samotij~\cite{BMS} and of Saxton and Thomason~\cite{ST} provides a useful framework for controlling independent sets in uniform hypergraphs. The method has been widely used in Ramsey theory and additive combinatorics. We shall use the following one-step container lemma, whose proof can be found in \cite[Lemma 2.3]{QL}.
\begin{lemma}[Hypergraph containers \cite{BMS,ST}]\label{ct:containers}
For every integer $m\ge 2$ there are constants $c_{m},C_{m},d_{m},\tau_m>0$ with the following property. Let $G$ be an $m$-graph with average degree $d>d_{m}$, and let $0<\tau\le\tau_m$.
If
\[\Delta_{j}(G)\le c_{m}d\tau^{j-1}\qquad (2\le j \le m),\]
then there is a family $\mathcal{C}$ of subsets of $V(G)$ such that
\begin{enumerate}
    \item [(i)] every independent set in $G$ is contained in some $C \in \mathcal{C}$;
    \item[(ii)] $e(G[C])\le e(G)/2 $ for every $C \in \mathcal{C}$;
    \item[(iii)]$\log|\mathcal C|\le C_{m}|V(G)|\tau \log(1/\tau)$.
\end{enumerate}
\end{lemma}

\section{The Construction}\label{sec:construction}
The following is the key lemma for proving Theorems~\ref{thm:half} and \ref{thm:localbase}, and its proof will be given in a subsequent section.
\begin{lemma}\label{lem:general}
Fix integers $k\ge3$ and $2\le r<k$. Put $m=\binom kr$ and $\beta=1/(48m^2)$. For every sufficiently large $n$, with $q=\lfloor n^{1/(16m)}\rfloor$, there is an oriented palette $\mathcal P\subseteq[q]^m$ whose support hypergraph is linear and Berge-triangle-free, and which meets every product $A_1\times\cdots\times A_m$ with $|A_i|\ge\lceil q^{1-1/(12m)}\rceil$ for all $i$. For every \(N\ge n\), let \(f:\binom{[N]}{r}\to[q]\) be chosen uniformly at random.
Then, for every fixed \(W\in\binom{[N]}{n}\),
\begin{equation}\label{eq:general-probability}
 \mathbb P\bigl(e(R(f,\mathcal P)[W])=0\bigr)
 \le q^{-\beta\binom nr}.
\end{equation}
Consequently, for a fixed family $\mathcal F$ of $k$-graphs, if the $k$-graph $R(f,\mathcal P)$ is $\mathcal F$-free for every $N$, every labeling $f$, and every oriented palette whose support hypergraph is linear and Berge-triangle-free, then, for some constant $c_{k,r}>0$ and all sufficiently large $n$,
\begin{equation}\label{eq:general-Ramsey}
 r_k(\mathcal F,n)>2^{c_{k,r}n^{r-1}\log n}.
\end{equation}
\end{lemma}

\subsection{Allowed patterns}
We give the palette construction for an arbitrary number $m$ of coordinates. The random selection and deletion argument extends the construction for triples in Conlon, Fox, Gunby, He, Mubayi,~Suk and Verstra\"ete~\cite[Lemma~2.4]{CFGHMSV}. Its intersection conditions are used for deterministic avoidance and its product condition is used for counting.

\begin{lemma}\label{lem:patterns}
Fix $m\ge3$ and put $\eta=1/(12m)$. For every sufficiently large integer $q$, with $L=\lceil q^{1-\eta}\rceil$, there is a family $\mathcal P\subseteq[q]^m$ with the following properties.
\begin{enumerate}
\renewcommand{\labelenumi}{(\roman{enumi})}
\renewcommand{\theenumi}{\roman{enumi}}
\item\label{it:unique} Each vector has distinct coordinates, and no two distinct vectors have the same support.
\item\label{it:linear} The support hypergraph $\mathcal A(\mathcal P)$ is linear and has no Berge triangle.
\item\label{it:product} Whenever \(A_1,\ldots,A_m\subseteq[q]\) satisfy \(|A_i|\ge L\) for all \(i\) (the sets \(A_i\) need not be disjoint), the product \(A_1\times\cdots\times A_m\) meets \(\mathcal P\).
\end{enumerate}
\end{lemma}
\begin{proof}
Begin with the color set $[2q]$. Independently select each $m$-subset with probability $p=q^{-(m-1)+1/6}$. For each selected set choose one of its $m!$ orderings uniformly and independently. Let $\mathcal P_0$ be the resulting vectors and let $\mathcal A_0$ be their support hypergraph.

There are $O_m(q^{2m-2})$ pairs of distinct $m$-sets meeting in at least two vertices. The expected number of selected such pairs is $O_m(q^{2m-2}p^2)=O_m(q^{1/3})$. A Berge triangle can be specified by its three distinct intersection vertices and at most $m-2$ additional vertices in each of its three distinct edges. There are $O_m(q^{3m-3})$ such specifications, each selected with probability $p^3$. Their expected number is $O_m(q^{3m-3}p^3)=O_m(q^{1/2})$. If $Z$ counts all these pairs and triangle specifications, then $\mathbb EZ=O_m(q^{1/2})$, and Markov's inequality gives $\mathbb P(Z\ge q)\le O_m(q^{-1/2})=o(1)$.

Fix sets $A_1,\ldots,A_m\subseteq[2q]$ of size $L$. Among their $L^m$ ordered tuples, at most $\binom m2L^{m-1}$ have repeated coordinates. For sufficiently large $q$, at least $L^m-\binom m2L^{m-1}\ge L^m/2$ therefore have distinct coordinates. For each $B\in\binom{[2q]}m$, let $t_B$ count its orderings in this product, and so $\sum_Bt_B\ge L^m/2$. On the other hand, for a fixed support \(B\), there are two mutually exclusive ways in which \(B\) can fail to contribute a hit. First, \(B\) may not be selected at all. This happens with probability \(1-p\), and then \(B\) certainly contributes no hit. Second, \(B\) may be selected, which happens with probability \(p\), but its uniformly random ordering may fail to lie in the product \(A_1\times\cdots\times A_m\). Since exactly \(t_B\) of the \(m!\) orderings of \(B\) lie in the product, the conditional probability of failure given that \(B\) is selected is \(1-\frac{t_B}{m!}\). Independence across supports gives
\begin{align}\label{eq:box-miss}
 \mathbb P\bigl(\mathcal P_0\cap(A_1\times\cdots\times A_m)=\varnothing\bigr)
 &\le \prod_{B\in\binom{[2q]}m}\left((1-p)+p\left(1-\frac{t_B}{m!}\right)\right)\notag=\prod_{B\in\binom{[2q]}m}\left(1-\frac{pt_B}{m!}\right)\notag\\
 &\le\exp\left(-\frac p{m!}\sum_Bt_B\right)
 \le\exp\left(-\frac{pL^m}{2m!}\right).
\end{align}
Since $m\eta=1/12$, we have $pL^m\ge q^{13/12}$. There are at most $\binom{2q}L^m\le(2q)^{mL}$ choices of the sets. The union bound shows that the probability that some product is missed is at most $\exp(mL\ln(2q)-q^{13/12}/(2m!))=o(1)$.

Choose an outcome with $Z<q$ in which every such product is met. From the union of the edges in each counted pair or triangle choose one vertex, and let $W$ be the set of chosen vertices. Then $|W|<q$. Delete all vertices in $W$ and all edges meeting $W$. Each counted configuration loses an edge, so the remaining hypergraph is linear and Berge-triangle-free. Choose any $q$-element set $Q\subseteq[2q]\setminus W$ and retain exactly the vectors whose supports lie in $Q$.

Every product of $L$-subsets of $Q$ has an original allowed vector, and that vector is retained because all its coordinates lie in $Q$. Larger products contain such a product. Finally, relabel $Q$ as $[q]$. Distinct coordinates and unique orientations are preserved, completing the proof.
\end{proof}

\subsection{Common supports and vertex positions}
We first record the consequence of assigning a unique orientation to each palette support. 

\begin{lemma}\label{str:positions}
Let $2\le r<k$, let $m=\binom{k}{r}$, and suppose that every vector in $\mathcal P\subseteq[q]^m$ has distinct coordinates and that different vectors have different supports. If two edges $E,F$ of $R(f,\mathcal P)$ have the same palette support, then every vertex $v$ that belongs to a common $r$-face of $E$ and $F$ occupies the same position in $E$ and $F$. In particular, this conclusion holds for every $v\in E\cap F$ whenever $|E\cap F|\ge r$.
\end{lemma}

\begin{proof}
Let \(E\) and \(F\) be two edges in \(R(f,\mathcal P)\) such that
\[
\{f(\sigma_1(E)),\ldots,f(\sigma_m(E))\}
=
\{f(\sigma_1(F)),\ldots,f(\sigma_m(F))\}.
\]
Then
\[
\bigl(f(\sigma_1(E)),\ldots,f(\sigma_m(E))\bigr)
=
\bigl(f(\sigma_1(F)),\ldots,f(\sigma_m(F))\bigr),
\]
since different vectors in \(\mathcal P\) have different supports. Let \(\sigma\subseteq E\cap F\) be an \(r\)-face. Then \(\sigma=\sigma_i(E)=\sigma_j(F)\) for some \(i\) and \(j\). Moreover, \(i=j\), since every vector in \(\mathcal P\subseteq[q]^m\) has distinct coordinates. For any \(v\in\sigma\), the vertices of \(\sigma\) have the same relative order in the global vertex ordering, and hence each of them occupies the same position in \(E\) and in \(F\). If \(|E\cap F|\ge r\), then every vertex of the intersection belongs to such an \(r\)-face.
\end{proof}

\section{Excluding Forbidden Configurations}\label{sec:avoidance}

Throughout this section, \bhz{let $m=\binom{k}{r}$ and $\mathcal P\subseteq[q]^m$ be any palette satisfying} properties~\textup{(\ref{it:unique})} and~\textup{(\ref{it:linear})} of Lemma~\ref{lem:patterns}. \bhz{These two properties suffice for all the exclusion results in this section, while property~\textup{(\ref{it:product})} is needed only for the subsequent counting argument.}

\subsection{Three-edge configurations}

We begin with the configurations $\Delta_{k,a}$. The proof separates the case in which two edges share several $r$-faces from the case in which they share exactly one. Recall that we call all edges of \(R(f,\mathcal P)\) red and all other edges blue.

\begin{lemma}\label{str:delta}
Let $k\ge3$, $2\le r<k$, and $1\le a\le\min\{k-r,\lfloor k/2\rfloor\}$. Then $R(f,\mathcal P)$ is $\Delta_{k,a}$-free for every labeling $f:\binom{[N]}r\to[q]$.
\end{lemma}

\begin{proof}
Suppose that $R(f,\mathcal P)$ contains three edges $E_1=U\cup X\cup Y$, $E_2=U\cup X\cup Z$, and $E_3=U\cup Y\cup Z$, where $U,X,Y,Z$ are pairwise disjoint, $|U|=k-2a$, and $|X|=|Y|=|Z|=a$. Each pair of edges meets in $k-a\ge r$ vertices. Let $A_i=\{f(\sigma_1(E_i)),\ldots,f(\sigma_m(E_i))\}$ be the palette support of $E_i$. 

If \(k-a>r\), then \(|\sigma(E_i)\cap \sigma(E_j)|\ge \binom{k-a}{r}\ge 2\), that is, for every distinct \(i,j\in[3]\), the sets \(E_i\) and \(E_j\) share at least two distinct \(r\)-faces. Suppose that \(\sigma_s(E_i),\sigma_t(E_i)\in \sigma(E_i)\cap \sigma(E_j)\). Then \(f(\sigma_s(E_i))\neq f(\sigma_t(E_i))\), since both are contained in a red edge whose color vector has pairwise distinct coordinates by property~(\ref{it:unique}). Hence \(|A_i\cap A_j|\ge 2\). By the linearity property of the support hypergraph, property~(\ref{it:linear}), we obtain \(A_1=A_2=A_3\).

Suppose instead that $k-a=r$. Put $\alpha=f(U\cup X)$, $\beta=f(U\cup Y)$, and $\gamma=f(U\cup Z)$. Any two of these three distinct faces lie together in one of the red edges, so $\alpha,\beta,\gamma$ are pairwise distinct. Moreover, $\{\alpha,\beta\}\subseteq A_1$, $\{\alpha,\gamma\}\subseteq A_2$, and $\{\beta,\gamma\}\subseteq A_3$. If the supports were pairwise distinct, they would form a Berge triangle with intersection vertices $\alpha,\beta,\gamma$. If exactly two supports were equal, the third would meet that common support in at least two of these colors, contrary to linearity. Hence $A_1=A_2=A_3$ in this case as well.

Lemma~\ref{str:positions} now shows that every vertex has the same position in every edge containing it. For \(Q\in\{U,X,Y,Z\}\), let \(P_Q(E_i)\) denote the set of positions of \(Q\) in the edge \(E_i\). This is well-defined even when \(U\) is empty. The descriptions of \(E_1\) and \(E_2\) yield $P_U(E_1)\cup P_X(E_1)=P_U(E_2)\cup P_X(E_2)$.
Then
\[
P_Y(E_1)=[k]\setminus(P_U(E_1)\cup P_X(E_1))
=[k]\setminus(P_U(E_2)\cup P_X(E_2))
=P_Z(E_2).
\]
Similarly, \(P_Y(E_1)=P_Y(E_3)\) and \(P_Z(E_2)=P_Z(E_3)\). Consequently, \(P_Y(E_3)=P_Z(E_3)\). However, since \(Y\) and \(Z\) are nonempty disjoint subsets of \(E_3\), the sets \(P_Y(E_3)\) and \(P_Z(E_3)\) are nonempty and disjoint. This contradiction proves the lemma.
\end{proof}

The next lemma includes the pair-labeling case $r=2$. It is the step that connects the simultaneous avoidance construction to half-graphs.

\begin{lemma}\label{str:overlap}
Let $k\ge4$ and $2\le r\le k-2$. Then $R(f,\mathcal P)$ is $T_{k,r}$-free for every labeling $f:\binom{[N]}r\to[q]$.
\end{lemma}

\begin{proof}
Suppose that the three edges $A=S\cup\{x\}$, $B=S\cup\{y\}$, and $C=U\cup\{x,y\}\cup W$ are red, where $|S|=k-1$, $U\subseteq S$ has size $r-1$, and $W$ is disjoint from $S\cup\{x,y\}$ and has size $k-r-1$.

The edges $A$ and $B$ share $\binom{k-1}{r}\ge2$ distinct $r$-faces. It follows from the first part of property~(\ref{it:linear}) that
$\{f(\sigma_1(A)),\ldots,f(\sigma_m(A))\}=\{f(\sigma_1(B)),\ldots,f(\sigma_m(B))\}$,
which, together with property~(\ref{it:unique}), implies that
\[
(f(\sigma_1(A)),\ldots,f(\sigma_m(A)))
=
(f(\sigma_1(B)),\ldots,f(\sigma_m(B))).
\]
Furthermore, Lemma~\ref{str:positions} implies that every vertex of $S$ has the same position in $A$ and $B$. The vertices $x$ and $y$ must consequently occupy the same remaining position. Thus, the $r$-faces $U\cup\{x\}$ and $U\cup\{y\}$ occur at the same coordinate in the equal color vectors of $A$ and $B$. Specifically, $\sigma_i(A)=U\cup\{x\}$ and $\sigma_i(B)=U\cup\{y\}$ for some $i\in[m]$. Hence $f(U\cup\{x\})=f(U\cup\{y\})$. These are two distinct $r$-faces of the red edge $C$, contradicting the requirement that its vector have distinct coordinates.
\end{proof}

Combining Lemma~\ref{str:delta} and Lemma~\ref{str:overlap}, we directly obtain the following corollary.
\begin{corollary}\label{str:family}
For every $k\ge3$ and $2\le r<k$, the hypergraph $R(f,\mathcal P)$ is $\mathcal F_{k,r}$-free for every labeling $f:\binom{[N]}r\to[q]$.
\end{corollary}

As mentioned in the introduction, the two special elements in $\mathcal F_{k,r}$ will be used to prove Theorem~\ref{thm:half} and Theorem~\ref{thm:localbase}, respectively.
 More specifically, for $k\ge4$, the half-graph $B^{(k)}$ contains a copy of $T_{k,2}$. Moreover, $B^{(3)}\cong\Delta_{3,1}$, and every three distinct $k$-edges on a set of $k+1$ vertices form a copy of $\Delta_{k,1}$. Consequently, a pair-labeling construction satisfying the hypotheses of this section is $B^{(k)}$-free for every $k\ge3$, and a $(k-1)$-face-labeling construction has at most two red edges on every $(k+1)$-set.

\subsection{Connectivity of intersections}

The next proposition gives a sufficient condition for deterministic exclusion in terms of the intersection graphs defined in Section~\ref{sec:introduction}. It uses linearity of the palette support hypergraph, but does not require the absence of Berge triangles.
For an ordered $k$-set $E$ and $v\in E$, write $\rho_E(v)=1+|\{w\in E:w<v\}|$ for the position of $v$ in $E$. 

Recall that, for a $k$-graph $F$ and an integer $s\in[k]$, we let $\Gamma_s(F)$ denote the graph with vertex set $E(F)$ in which two distinct edges of $F$ are adjacent when they meet in at least $s$ vertices. Let $\Gamma_s(F;v)$ be its induced subgraph on the edges containing $v$. A graph is connected if any two of its vertices are joined by a path. In particular, a one-vertex graph is connected. For $k\ge4$ and $2\le r\le k-2$, let $\mathcal H_{k,r}$ be the family of non-$k$-partite $k$-graphs $F$ with no isolated vertices such that $\Gamma_{r+1}(F)$ is connected and $\Gamma_r(F;v)$ is connected for every $v\in V(F)$.
\begin{proposition}\label{oc:criterion}
Fix $k\ge4$ and $2\le r\le k-2$. Let $F$ be a $k$-graph without isolated vertices such that $\Gamma_{r+1}(F)$ is connected and $\Gamma_r(F;v)$ is connected for every $v\in V(F)$. If $\mathcal P$ is an oriented $r$-face palette with linear support hypergraph, then every copy of $F$ in $R(f,\mathcal P)$ induces a $k$-partition of $F$. Consequently, $R(f,\mathcal P)$ is $\mathcal H_{k,r}$-free for every labeling $f$.
\end{proposition}

\begin{proof}
Suppose that $F$ appears as a subhypergraph. If $E,E'\in E(F)$ have at least $r+1$ common vertices, then they share at least $\binom{r+1}{r}=r+1$ distinct $r$-faces. The labels of these faces are distinct because they occur in a red edge. Hence $$\{f(\sigma_1(E)),\ldots,f(\sigma_m(E))\}=\{f(\sigma_1(E')),\ldots,f(\sigma_m(E'))\},$$ since $\mathcal P$ is an oriented $r$-face palette with linear support hypergraph. Connectivity of $\Gamma_{r+1}(F)$ now implies that all edges of $F$ have the same palette support.

Fix $v\in V(F)$. If two edges containing $v$ meet in at least $r$ vertices, Lemma~\ref{str:positions} shows that $v$ has the same position in both. Connectivity of $\Gamma_r(F;v)$ therefore makes $\rho_E(v)$ independent of the choice of the edge $E$ containing $v$. Since $F$ has no isolated vertices, we may define $c(v)=\rho_E(v)\in[k]$ for every $v\in V(F)$. For every edge $E=\langle v_1,\ldots,v_k\rangle$, we have $c(v_i)=i$ for $i\in[k]$. The classes of $c$ form a $k$-partition of $F$. This excludes every member of $\mathcal H_{k,r}$.
\end{proof}

\section{Counting Independent Sets}\label{sec:counting}
In this section, we will use the container method to count independent sets of size $n$ in $R(f,\mathcal P)$.
\subsection{The auxiliary hypergraph $\mathcal G$}
Fix a set \(W\subset [N]\) with \(|W|=n\). 
For the remainder of the proof, let
\begin{equation}\label{canshu}
    m=\binom{k}{r},~~\eta=\frac{1}{12m}~~\text {and}~~\beta=\frac{\eta}{4m}=\frac{1}{48m^2}.
\end{equation}

Given a sufficiently large integer \(n\), set
\(
q=\left\lfloor n^{1/(16m)}\right\rfloor\) and \(L=\left\lceil q^{1-\eta}\right\rceil,
\)
and choose a family \(\mathcal P\) as in Lemma~\ref{lem:patterns}.

Define an auxiliary $m$-graph $\mathcal G=\mathcal G(n,q)$ on $\binom{W}{r}\times[q]$.  For every $E \in \binom{W}{k}$ and every $(a_1,...,a_m)\in \mathcal P$, include the edge
\begin{equation}\label{eq:auxiliary-edge}
 \{(\sigma_i(E),a_i):i \in [m]\}.
\end{equation}
For a local coloring \(g:\binom{W}{r}\to[q]\) and \(\mathcal{P}\subseteq [q]^m\), define a \(k\)-graph \(R(g,\mathcal{P})\) on \(W\) by declaring that \(E'=\langle v_1,\ldots,v_k\rangle\) is an edge if and only if \(\bigl(g(\sigma_1(E')),\ldots,g(\sigma_m(E'))\bigr)\in\mathcal P\). In particular, if \(W=\{w_1,\ldots,w_n\}\) with \(w_1<\cdots<w_n\) and \(g(\sigma)=f(\sigma)\) for all \(\sigma\in \binom{W}{r}\), then \(R(g,\mathcal{P})=R(f,\mathcal{P})[W]\), the induced subhypergraph of \(R(f,\mathcal{P})\) on the vertex set \(W\).

Let $$I_g=\left\{(\sigma,g(\sigma)):\sigma\in \binom{W}{r}\right\}.$$

A local coloring $g$ corresponds to the set $I_g=\{(e,g(e)):e\in\binom{W}r\}$. By the definition of $\mathcal G$, we have the following lemma.

\begin{lemma}\label{independent set}
The set $I_g$ is independent in $\mathcal{G}$ if and only if  $R(g,\mathcal{P})$ is an edgeless hypergraph.
\end{lemma}
\begin{proof}
First suppose that $R(g,\mathcal P)$ contains an edge $E'=\langle v_1,\ldots,v_k\rangle$. Then
\begin{equation*}
\bigl(g(\sigma_1(E')),\ldots,g(\sigma_m(E'))\bigr)\in\mathcal P.
\end{equation*}
Hence $\{(\sigma_i(E'),g(\sigma_i(E'))):i\in[m]\}$ is an edge of $\mathcal G$ contained in $I_g$, so $I_g$ is not independent.

Conversely, suppose that $I_g$ is not independent in $\mathcal G$. Then it contains an edge of $\mathcal G$ arising from some $E\in\binom{W}{k}$ and $\mathbf a=(a_1,\ldots,a_m)\in\mathcal P$. Since this edge is contained in $I_g$, we have $a_i=g(\sigma_i(E))$ for every $i\in[m]$. Thus $E$ is an edge of $R(g,\mathcal P)$.
\end{proof}

\subsection{Counting local colorings with no edges}
Recall that the family \(\mathcal{P}\) is chosen as in Lemma~\ref{lem:patterns}.
\begin{lemma}\label{lem:probability}
For all sufficiently large \(n\), the number of local colorings \(g:\binom{W}{r}\to[q]\) such that \(R(g,\mathcal{P})\) is edgeless is at most \(q^{(1-\beta)\binom{n}{r}}\).
\end{lemma}

\begin{proof}
Set $\rho=\frac{1}{2m}$ and $D=\binom nr$. By Lemma~\ref{lem:supersaturation}, applied with $\varepsilon=(1-\frac{r}{2k})/{\binom{k-1}{r-1}}>0$, there is a constant $\kappa =\kappa (k,r)>0$ such that every $r$-graph $J$ on $n$ vertices with $e(J)\ge (1-\rho)D$ contains at least $\kappa n^k$ copies of $K_k^{(r)}$. Set $\zeta=\frac{\kappa}{2}.$

\medskip\noindent\textbf{Container applicability.}
We shall apply Lemma~\ref{ct:containers} recursively to induced subhypergraphs of $\mathcal{G}$. Let $C \subseteq V(\mathcal{G})$ satisfy $|E(\mathcal{G}[C])|\geq \zeta n^{k}$, and put $G:= \mathcal {G}[C]$. Let $d(G)$ denote the average degree of $G$. We have $|V(G)|\le |V(\mathcal G)|=qD$, therefore the average degree $d$ of $G$ satisfies
\begin{equation}\label{eq:average-degree}
 d=d(G)=\frac{m|E(G)|}{|V(G)|}\ge \frac{m\zeta n^k}{qD}\ge a_{k,\zeta}\frac{n^{k-r}}q.
\end{equation}
As $q=\lfloor n^{1/(16m)}\rfloor$, this implies $d(G)\ge a_{k,\zeta}n^{k-r-1/16m}\ge a_{k,\zeta}n^{1-1/16m} \to \infty$. Thus $d(G)>d_{m}$ in Lemma \ref{ct:containers} holds for all sufficiently large $n$.

\begin{claim}\label{cl:j-degree}
For every $2\le j\le m$,
\begin{equation}\label{eq:j-degree}
\Delta_{j}(G)\le n^{k-r-1}q^{m-j}.
\end{equation}
\end{claim}

\begin{proof}

Fix a \(j\)-set \(S\subseteq V(\mathcal{G})\) with \(2\le j\le m\). Each vertex of \(\mathcal{G}\) has the form \((e,c)\), where \(e\in \binom{W}{r}\) and \(c\in [q]\). Suppose $S$ contains two vertices $(e,c)$ and $(e,c')$ with $c\ne c'$. By the definition of $\mathcal G$, an edge contains at most one auxiliary vertex with any given first coordinate $e$. Hence no edge of $\mathcal G$ can contain both $(e,c)$ and $(e,c')$, and therefore no edge can contain $S$. Hence the codegree of $S$ is zero. We may therefore assume that the vertices of \(S\) involve \(j\) distinct $r$-sets from \(W\). Let \(X\subseteq W\) be the set of all elements of these $r$-sets. If $|X|>k$, the codegree of $S$ is zero. Hence we may assume $|X|\le k$. Since \(j\ge 2\) and the sets are distinct, we have \(|X|\ge r+1\). Every edge of \(\mathcal{G}\) containing \(S\) comes from some \(k\)-set \(T\subseteq W\). This set \(T\) must contain \(X\). Therefore, the number of possible choices for \(T\) is at most \(\binom{n-|X|}{k-|X|}\le n^{k-|X|}\le n^{k-r-1}\). For each such \(k\)-set, the prescribed pairs fix \(j\) color coordinates. The remaining \(m-j\) coordinates have at most \(q^{m-j}\) completions. Taking the maximum over all $j$-sets $S$, we obtain
\[
\Delta_j(\mathcal G)\le n^{k-r-1}q^{m-j}.
\]

As \(\Delta_{j}(G)\le \Delta_j(\mathcal{G})\), this proves~\eqref{eq:j-degree}. If \(S\) is not contained in any edge, then \(S\) has codegree zero and satisfies the same bound.
\end{proof}

Set $\tau=n^{-1/(4m)}$. Since $\tau\to 0$ as $n\to \infty$, for sufficiently large $n$, $\tau$ satisfies $\tau\le\tau_m$. Equations~\eqref{eq:average-degree} and~\eqref{eq:j-degree} give, for every $2\le j\le m$,
\begin{equation}\label{eq:codegree-ratio}
 \frac{\Delta_j(G)}{d\tau^{j-1}}
 \le a_{k,\zeta}^{-1}n^{-1}q^{m-j+1}\tau^{-(j-1)}
 \le a_{k,\zeta}^{-1}n^{-1+(m+3j-3)/(16m)}\le a_{k,\zeta}^{-1}n^{-3/4-3/(16m)}.
\end{equation}
We obtain that 
\[\max\limits_{2\leq j\leq m}\frac{\Delta_j(G)}{d(G)\tau^{j-1}}=O_k(n^{-3/4-3/(16m)}) = o(1).\]
All implicit constants here depend only on the fixed parameter $k$. Hence, for all sufficiently large $n$, Lemma~\ref{ct:containers} applies to every induced subhypergraph $\mathcal{G}[C]$ with $e(\mathcal{G}[C])\ge\zeta n^k$.

\medskip\noindent\textbf{Number of terminal containers.}
Starting from $V(\mathcal{G})$, apply Lemma~\ref{ct:containers} to every current container $C$ until $e(\mathcal{G}[C])<\zeta n^{k}$. Since 
\[e(\mathcal{G})=\binom{n}{k}|\mathcal{P}|,\qquad |\mathcal{P}|\le q^{m}\] and the number of edges is divided by at least two at each refinement step, the number of steps before reaching any terminal container is $O(\log{q})$.

At each step of the recursion, the covering property is preserved. To see this, let $I$ be an independent set of $\mathcal{G}$ contained in a current container $C$. Since $I$ is also independent in the induced subhypergraph $\mathcal{G}[C]$, Lemma~\ref{ct:containers} guarantees that $I$ is contained in one of the subcontainers generated from $C$. Repeating this argument recursively, we obtain a family $\mathcal{C}$ of terminal containers such that every independent set of $\mathcal{G}$ is contained in some $C\in \mathcal{C}$, and 
\[e(\mathcal{G}[C])<\zeta n^k\qquad \text{for every }C \in \mathcal{C}.\]
At each application of Lemma~\ref{ct:containers}, the logarithm of the number of new containers is at most
\[C_{m}|V(\mathcal{G})|\tau \log{(1/\tau)}.\]
As there are at most $O(\log{q})$ steps, we have
\begin{equation}\label{eq:container-asymptotic}
    \log{|\mathcal{C}|}\le O(Dq\log {q}\tau\log{(1/\tau)})=O(n^{r-3/(16m)}(\log{n})^{2})=o(D\log{q}).
\end{equation}
Here $q=\left\lfloor n^{\frac{1}{16m}}\right\rfloor$, $\tau=n^{-\frac{1}{4m}}$ and $D=\binom{n}{r}$, while $\log q=\Theta_k(\log n)$. The ratio in the last comparison is $O_k(n^{-\frac{3}{16m}}\log n)$, which tends to zero. Thus, $|\mathcal{C}|=q^{o(D)}$.

\medskip\noindent\textbf{Terminal containers have many small lists.}
Fix a terminal container $C\in\mathcal C$. For each original $r$-set  $\sigma\in\binom{W}r$, define its color list $A_{\sigma}(C)=\{a\in[q]:(\sigma,a)\in C\}$. Call $\sigma$ large if $|A_{\sigma}(C)|\geq L$. Let $J_{C}$ be the $r$-graph on $W$ whose edges are precisely the large sets.

\begin{claim}\label{cl:supersaturation-application}
For every terminal container $C\in \mathcal{C}$,
\[e(J_{C})<(1-\rho)D.\]
\end{claim}

\begin{proof}
Suppose, for a contradiction, that $e(J_{C})\geq (1-\rho)D$. By the choice of $\kappa$, the $r$-graph $J_{C}$ contains at least $\kappa n^k$ distinct copies of $K_{k}^{(r)}$. Let $v_1<\cdots<v_k$ be the vertices of one such copy. Since every $r$-set $\sigma_{i}$ is large, for each $i\in [m]$, all lists $A_{\sigma_{i}}(C)$ have size at least $L$. The property~\textup{(\ref{it:product})} of Lemma~\ref{lem:patterns} gives $\boldsymbol{a}=(a_{i})_{i\in [m]} \in \mathcal{P}\cap (A_{\sigma_{1}}(C)\times\cdots \times A_{\sigma_{m}}(C))$. Hence it gives an edge $\left\{(\sigma_{i},a_{i}):i\in [m]\right\}$
of $\mathcal{G}[C].$
Distinct $k$-vertex sets spanning copies of $K_{k}^{(r)}$ in $J_{C}$ give distinct edges of $\mathcal{G}[C]$. Therefore
\[e(\mathcal{G}[C])\ge \kappa n^{k}> \zeta n^{k},\]
contradicting the terminal property of $C$. This proves the claim.
\end{proof}

\medskip\noindent\textbf{Count bad local colorings.} Fix a terminal container \(C \in \mathcal{C}\). Recall that
\[
I_g = \left\{ (\sigma, g(\sigma)) : \sigma \in \binom{W}{r} \right\}, \qquad A_{\sigma}(C) = \{ a \in [q] : (\sigma, a) \in C \}.
\]

For each $r$-set \(\sigma \in \binom{W}{r}\), the set \(I_g\) contains the single vertex \((\sigma, g(\sigma))\) of \(V(\mathcal{G})\) lying over \(\sigma\). Therefore \(I_g \subseteq C\) if and only if \(g(\sigma) \in A_{\sigma}(C)\) for every \(\sigma \in \binom{W}{r}\). Hence the number of colorings \(g\) with \(I_g \subseteq C\) is at most \(\prod_{\sigma \in \binom{W}{r}} |A_{\sigma}(C)|\).

Every list has size at most $q$, while a list of size less than $L$ has size at most $L-1<q^{1-\eta}$. Thus, the number of local colorings $g$ with $I_{g}\subseteq C$ is
\begin{equation}\label{eq:coloring-count}
 \prod_{\sigma\in\binom{W}r}|A_{\sigma}(C)|\le q^{(1-\rho)D}(q^{1-\eta})^{\rho D}=q^{(1-\rho)D+\rho (1-\eta)D}
 = q^{(1-\eta/(2m))D}.
\end{equation}
This inequality also holds if one of the lists is empty, when the product is zero.

Thus, for each fixed terminal container \(C\), the number of local colorings \(g\) with \(I_g \subseteq C\) is at most \(q^{(1-\eta/(2m))D}\). By Lemma \ref{independent set}, every local coloring \(g\) for which \(R(g,\mathcal{P})\) contains no edge gives an independent set \(I_g\) in \(\mathcal{G}\). Since every independent set is contained in some terminal container, the total number of such local colorings is at most
\[
|\mathcal{C}|q^{(1-\eta/(2m))D} = q^{o(D)}q^{(1-\eta/(2m))D}.
\]
By~\eqref{eq:container-asymptotic}, for sufficiently large $n$ we have $|\mathcal C|\le q^{\eta D/(4m)}$. Hence the total number of local colorings \(g\) for which \(R(g,\mathcal{P})\) contains no edge is at most
\[
q^{\eta D/(4m)}q^{(1-\eta/(2m))D} =q^{(1-\eta/(4m))D}\overset{\eqref{canshu}}{=}q^{(1-\beta)D}.
\]
This proves the lemma.
\end{proof}

\section{Proofs of the Main Results}\label{sec:main-proof}
\subsection{Applications of Oriented Palette Constructions}
We first prove Lemma~\ref{lem:general}.
\begin{proof}[Proof of Lemma~\ref{lem:general}]
Lemma~\ref{lem:patterns} supplies the asserted palette. On any fixed $n$-set $W$, put $g=f|_{\binom{W}{r}}$. By Lemma~\ref{lem:probability},
\begin{align*}
 \mathbb P\bigl(e(R(f,\mathcal P)[W])=0\bigr)
 &=\mathbb P\bigl(e(R(g,\mathcal P))=0\bigr)\\
 &\le\frac{q^{(1-\beta)\binom nr}}{q^{\binom nr}}
 =q^{-\beta\binom nr}.
\end{align*}
Here all $q^{\binom nr}$ local colorings on $\binom{W}{r}$ are equally likely. This gives~\eqref{eq:general-probability} for every $N\ge n$.

Suppose now that $\mathcal F$ satisfies the hypothesis of Lemma~\ref{lem:general}. Put $D=\binom nr$ and take $N=\lfloor q^{\beta D/(2n)}\rfloor$. Let $Y$ count independent $n$-sets of $R(f,\mathcal P)$. Let $W$ be any fixed $n$-set. Then, by \eqref{eq:general-probability} and the linearity of expectation, it follows that
\begin{equation}\label{eq:first-moment}
 \mathbb EY\le\binom Nn q^{-\beta D}
 \le N^nq^{-\beta D}
 \le q^{-\beta D/2}<1.
\end{equation}
Thus there is some labeling $f$ such that $Y=0$. Its red hypergraph is $\mathcal F$-free by hypothesis, so $r_k(\mathcal F,n)>N$. For sufficiently large $n$, we have $D\ge n^r/(2r!)$ and $\log q\ge(\log n)/(32m)$. The floor in the definition of $N$ gives
\[
 \log N\ge\frac{\beta D}{2n}\log q-1
 \ge\frac{\beta}{128mr!}n^{r-1}\log n-1,
\] 
which proves~\eqref{eq:general-Ramsey}.
\end{proof}

We now prove Theorems~\ref{thm:family} and~\ref{oc:ramsey}.

\begin{proof}[Proof of Theorem~\ref{thm:family}]
By combining Corollary~\ref{str:family} with Lemma~\ref{lem:general} for \(\mathcal F=\mathcal F_{k,r}\), we obtain the desired result for every \(2\le r<k\).
\end{proof}

\begin{proof}[Proof of Theorem~\ref{oc:ramsey}]
For every palette in Lemma~\ref{lem:general}, Corollary~\ref{str:family} excludes every member of $\mathcal F_{k,r}$ and Proposition~\ref{oc:criterion} excludes every member of $\mathcal H_{k,r}$. Both assertions hold for the same labeling. Applying Lemma~\ref{lem:general} to their union gives the desired result.
\end{proof}

As mentioned in the introduction, Theorems~\ref{thm:half} and~\ref{thm:localbase} are direct consequences of Theorem~\ref{thm:family}.
\begin{proof}[Proof of Theorem~\ref{thm:half}]
For $k\ge4$, the half-graph $B^{(k)}$ contains a copy of $T_{k,2}$. The $\mathcal F_{k,2}$-free construction in Theorem~\ref{thm:family} is therefore $B^{(k)}$-free and gives $r_k(B^{(k)},n)>2^{c_kn\log n}$. For $k=3$, use $B^{(3)}=\Delta_{3,1}$ and Theorem~\ref{thm:family} with $r=2$.
\end{proof}

\begin{proof}[Proof of Theorem~\ref{thm:localbase}]
Since any three distinct $k$-edges on a $(k+1)$-set form a copy of $\Delta_{k,1}$, we have $r_k(k+1,3;n)=r_k(\Delta_{k,1},n)$. Corollary~\ref{cor:delta} with $a=1$ gives the lower bound $r_k(\Delta_{k,1},n)> 2^{\Omega_k(n^{k-2}\log n)}$, while the upper bound in~\eqref{r_k(k+1,3;n)} gives the matching upper bound.
\end{proof}

Recall that \(C_\ell^{(k)}\) is the tight cycle on vertices
\(v_0,\ldots,v_{\ell-1}\), whose edges are
\(E_i=\{v_i,\allowbreak v_{i+1},\allowbreak\ldots,\allowbreak v_{i+k-1}\}\)
for \(0\le i<\ell\), with subscripts interpreted modulo \(\ell\).
Recall also that \(F_{\ell,k,d}\) is obtained by deleting
\(E_0,\ldots,E_{d-1}\).

Next, we prove Corollary~\ref{oc:cycles}.
In the proof, subscripts are once again interpreted modulo \(\ell\).

\begin{proof}[Proof of Corollary~\ref{oc:cycles}]
Put \(r=k-d-1\), so that \(2\le r\le k-2\). It suffices to verify that
\(F_{\ell,k,d}\in\mathcal H_{k,r}\), since then
\[
r_k(F_{\ell,k,d},n)
\ge r_k(\mathcal H_{k,k-d-1},n)
\ge r_k(\mathcal F_{k,k-d-1}\cup\mathcal H_{k,k-d-1},n)
> 2^{c_{k,\ell,d}n^{k-d-2}\log n},
\]
where the last inequality follows from Theorem~\ref{oc:ramsey}. Taking $d=1$ yields the asserted bound for $C_\ell^{(k)}-e$.

For $d\le i<\ell-1$, the consecutive retained edges $E_i,E_{i+1}$ intersect in $k-1\ge r+1$ vertices. Hence $E_d,E_{d+1},\ldots,E_{\ell-1}$ is a path containing every vertex of $\Gamma_{r+1}(F_{\ell,k,d})$, and this graph is connected.

Fix a vertex \(v_j\). In \(C_\ell^{(k)}\), we have \(v_j\in E_t\) for all
\(t\in[j-k+1,j]\), where \([j-k+1,j]\) is a cyclic interval of length \(k\)
(indices are understood modulo \(\ell\)). Note that
\[
[j-k+1,j]\cap[d,\ell-1]=
\begin{cases}
[\ell+j-k+1,\ell-1], & 0\le j\le d-1, \\[2mm]
[d,j]\cup[\ell+j-k+1,\ell-1], & j\in[d,k-2], \\[2mm]
[d,j], & j\in[k-1,d+k-2], \\[2mm]
[j-k+1,j], & j\in[d+k-1,\ell-1].
\end{cases}
\]
Consequently, the intersection of \([j-k+1,j]\) with the retained index
interval \([d,\ell-1]\) consists of either one interval or two intervals in
the usual linear order. Within each interval, consecutive edges intersect in
\(k-1\ge r+1>r\) vertices. If two intervals occur, then \(v_j\in E_d\) and
\(v_j\in E_{\ell-1}\). The assumption \(\ell\ge2k+1\) gives
$E_d\cap E_{\ell-1}=\{v_d,v_{d+1},\ldots,v_{k-2}\}$.
This intersection has \(k-d-1=r\) vertices. Hence \(E_d\) and
\(E_{\ell-1}\) are adjacent in \(\Gamma_r(F_{\ell,k,d};v_j)\). Therefore,
\(\Gamma_r(F_{\ell,k,d};v_j)\) is connected, and so
\(\Gamma_r(F_{\ell,k,d};v)\) is connected for every vertex \(v\). Finally,
there are no isolated vertices, since every vertex originally belongs to
\(k\) edges and only \(d<k\) edges have been deleted, so every vertex still
belongs to at least one retained edge.

It remains to show that \(F_{\ell,k,d}\) is not \(k\)-partite. Suppose, to the
contrary, that it is, and let \(c(v)\in[k]\) denote the color class of \(v\).
Define
\[
a_i=c(v_{i\bmod\ell})\qquad(d\le i\le \ell+k-2).
\]
For every \(d\le i\le \ell-1\), the retained edge \(E_i\) implies that
\(a_i,\ldots,a_{i+k-1}\) are pairwise distinct. Comparing the two consecutive
edges \(E_i\) and \(E_{i+1}\) gives, for \(d\le i\le \ell-2\),
\begin{equation}\label{oc:periodicity}
a_i=a_{i+k}.
\end{equation}
Hence the sequence \(\{a_i\}_{i\in[d,\ell+k-2]}\) is \(k\)-periodic, and each
period contains all \(k\) distinct classes. Note that
\(k-2,\ell+k-2\in[d,\ell+k-2]\), and
\[
a_{k-2}
=c(v_{(k-2)\bmod\ell})
=c(v_{(\ell+k-2)\bmod\ell})
=a_{\ell+k-2}.
\]
By \eqref{oc:periodicity}, distinct residues modulo \(k\) have distinct colors.
Therefore \(a_{k-2}=a_{\ell+k-2}\) forces
\[
k-2\equiv \ell+k-2\pmod{k},
\]
that is, \(k\mid \ell\), a contradiction. Therefore, \(F_{\ell,k,d}\) is not \(k\)-partite. This completes the proof.
\end{proof}

\subsection{The exact magnitude of $r_k(k+1,t;n)$ for $t\in[4,k-2]$}
In this subsection, we prove Lemma~\ref{su:oddstep} and Theorem~\ref{su:tower}.

\begin{proof}[Proof of Lemma~\ref{su:oddstep}]
Let \(N = r_{k-1}(k,3;n) - 1\) and $k=2h+1\ge7$. Suppose \(\varphi\) is a red/blue edge-coloring of \(K_N^{(k-1)}\) with at most two red edges on every $k$-set and no blue \(K_n^{(k-1)}\). Using \(\varphi\), we will construct a red/blue edge-coloring \(\chi\) for \(K_{2^N}^{(k)}\) that avoids a \((k+1)\)-vertex set spanning at least 4 red edges and avoids a blue \(K_{2kn}^{(k)}\). We take the vertex sets to be \(V(K_N^{(k-1)}) = [N]\) and \(V(K_{2^N}^{(k)}) = \{0, 1\}^N\).

We endow the vertex set \(V(K_{2^N}^{(k)})\) with the lexicographic order induced by their binary representations. Specifically, for two vertices \(a = (a(1), \dots, a(N))\) and \(b = (b(1), \dots, b(N))\), we have \(a < b\) if and only if at the first index \(i\) where they differ, \(a(i) = 0\) and \(b(i) = 1\) (with \(a(j) = b(j)\) for all \(1 \le j < i\)). For any distinct \(a, b \in V(K_{2^N}^{(k)})\), let \(\delta(a, b)\) denote this first differing index. We will utilize the following stepping-up properties (see \cite{D-H-L-W-4,G-R-S-1}).
\medskip

\textbf{Property I.} For every triple \(a < b < c\), \(\delta(a, b) \neq \delta(b, c)\).
\medskip

\textbf{Property II.} For \(a_1 < \dots < a_r\), \(\delta(a_1, a_r) = \min_{1 \le j \le r-1} \delta(a_j, a_{j+1})\), there is a unique $\delta_i$ which achieves the minimum.
\medskip

\textbf{Property III.} For $a_1 < \cdots < a_r$, suppose that $\delta_1,\ldots,\delta_{r-1}$ form a monotone sequence. Then for every subset of $s$ vertices $a_{i_1} < \cdots < a_{i_s}$, $\delta(a_{i_1},a_{i_2}),\delta(a_{i_2},a_{i_3}),\ldots,\delta(a_{i_{s-1}},a_{i_s})$ form a monotone sequence. Moreover, for every subset of $s-1$ such $\delta_j$'s, i.e. $\delta_{j_1},\delta_{j_2},\ldots,\delta_{j_{s-1}}$, there are $s$ vertices $a_{i_1},\ldots,a_{i_s}$ such that $\delta(a_{i_t},a_{i_{t+1}}) = \delta_{j_t}$.

\medskip

Before defining the coloring \(\chi\), we introduce some notation. Let \(V = \{0, 1\}^N\). Consider an \(s\)-element subset \(X = \{x_1, \dots, x_s\} \subset V\) with \(x_1 < \dots < x_s\), and define the sequence of integers \(\delta_i = \delta(x_i, x_{i+1})\) for \(i\in[s-1]\). We classify \(\delta_i\) based on its neighbors: it is a \textit{local minimum} if \(\delta_{i-1} > \delta_i < \delta_{i+1}\), a \textit{local maximum} if \(\delta_{i-1} < \delta_i > \delta_{i+1}\), and a \textit{local extremum} if it is either. Since adjacent terms in the sequence are always distinct (\(\delta_{i-1} \neq \delta_i\)), every nonmonotone sequence \(\{\delta_i\}_{i=1}^{s-1}\) must contain a local extremum. For shorthand, we write \(\delta(X) = \{\delta_i\}_{i=1}^{s-1}\). We say $\delta(X)$ is \textit{zigzag}, meaning $\delta_2, \delta_4,\ldots $ are local minima and $\delta_3, \delta_5,\ldots $ are local maxima. In other words, $\delta_1>\delta_2<\delta_3>\delta_4<\cdots$.

We are now ready to define the coloring \(\chi\) for the \(k\)-tuples in \(V = V(K_{2^N}^{(k)})\). Given an edge \(e = (x_1, \dots, x_k)\) where \(x_1 < \dots < x_k\), let \(\delta_i = \delta(x_i, x_{i+1})\). We use the notation $X\setminus\{x\}=X-x$. We assign \(\chi(e) = red\) if at least one of the following holds:
\begin{enumerate}
\item [($\Xi_1$)] The sequence $\delta_1,\ldots,\delta_{k-1}$ is monotone and $\varphi(\{\delta_1,\ldots,\delta_{k-1}\})=\text{red}$.
\item [($\Xi_2$)] $\delta(e-x_1)$ is zigzag and $\delta_1=\min_{i\in[k-1]}\delta_i$.
\item [($\Xi_3$)] $\delta(e-x_k)$ is zigzag and $\delta_{k-1}=\min_{i\in[k-1]}\delta_i$.
\end{enumerate}
All remaining edges are blue. The minimum in ($\Xi_2$) and ($\Xi_3$) is unique from \textbf{Property II}. Call an edge \(e\) satisfying \((\Xi_2)\) or \((\Xi_3)\) \emph{structural}, and an edge \(e\) satisfying \((\Xi_1)\) \emph{gap-monotone}. In particular, if \(\delta(e)\) is increasing, we call \(e\) \emph{gap-increasing}, and if \(\delta(e)\) is decreasing, we call \(e\) \emph{gap-decreasing}. 

For an ordered set $X=\{x_1<\cdots<x_s\}$ with $\delta(X)$, call $(x_i,x_{i+1})$ a \emph{pair} in $X$ if
\begin{equation}\label{eq:pair}
 \delta_i>\delta_{i-1}\ \text{when }i>1,
 \qquad \delta_i>\delta_{i+1}\ \text{when }i<s-1.
\end{equation}
For $i\in\{1,s-1\}$, it suffices to satisfy the sole existing neighboring inequality. Distinct pairs are disjoint, since two adjacent $\delta$ cannot both exceed each other. Call $X$ \emph{fully paired} if every point belongs to a pair.

\begin{fact}\label{three properties}
The following facts hold.
\begin{enumerate}
\item Any pair \((x_i,x_{i+1})\) contained in a subset \(X' \subseteq X\) remains a pair in \(X'\).
\item A \(2\ell\)-point set \(S\) is fully paired if and only if the sequence \(\delta(S)\) is zigzag, that is,  \(\delta_1>\delta_2<\delta_3>\cdots<\delta_{2\ell-1}\).
\end{enumerate}
\end{fact}

\begin{proof}[Proof of Fact \ref{three properties}]

For (1), if $x_j$ is the nearest retained point to the left of the pair $(x_i,x_{i+1})$, then
\[
\delta(x_j,x_i)=\min_{j\le a<i} \delta_a \le \delta_{i-1}<\delta_i,
\]
the right side is identical. For (2), disjoint consecutive pairs covering all points must be $(x_1,x_2),(x_3,x_4),\ldots,(x_{2\ell-1},x_{2\ell})$, and their defining inequalities are exactly those displayed.
\end{proof}
\medskip
\noindent\textbf{There is no blue $2kn$-set.}
Suppose for a contradiction that $S=\{a_1<\cdots<a_{2kn}\}$ is blue, and write $\delta_i=\delta(a_i,a_{i+1})$. Let $\delta_j$ be the unique minimum. The uniqueness of $\delta_{j}$ follows from Property II. 
Split $S$ into
\[
 S_L=\{a_1,\ldots,a_j\},\qquad
 S_R=\{a_{j+1},\ldots,a_{2kn}\}.
\]
Let $C$ be a larger side, so $q=|C|\ge kn$, and choose $z$ on the other side. By Property II, $\delta(z,c)=\delta_j$ for every $c\in C$.

There are at most $h-1$ pairs in $C$. Otherwise, we can choose $h$ pairs $(a_{i_1},a_{{i_1}+1}),\allowbreak\ldots,\allowbreak (a_{i_h},a_{{i_h}+1})$ in $C$, where $a_{i_1}<a_{{i_1}+1}<a_{i_2}<a_{{i_2}+1}<\cdots<a_{i_h}<a_{{i_h}+1}$. Let $A=\{a_{i_1}<a_{{i_1}+1}<a_{i_2}<a_{{i_2}+1}<\cdots<a_{i_h}<a_{{i_h}+1}\}$, then $|A|=2h=k-1$ since distinct pairs are disjoint. By Fact \ref{three properties} (1), all these pairs in $C$ remain pairs in $A$, so $A$ is fully paired. Fact \ref{three properties} (2) gives $\delta(A)$ is zigzag. If $C=S_L$, the set $A\cup\{z\}$ satisfies ($\Xi_3$). If $C=S_R$, it satisfies ($\Xi_2$). In both cases, this yields a red edge in \(S\), a contradiction.

Write \(C=\{c_1<\cdots<c_q\}\) and set \(\delta'_i=\delta(c_i,c_{i+1})\). Note that there is no integer \(\ell\) for which the sequence \(\{\delta'_i\}_{i=\ell}^{\ell+n-1}\) is monotone. Indeed, otherwise Property III would imply that \(\varphi\) colors every \((k-1)\)-tuple of the \(n\)-element set \(\{\delta'_i\}_{i=\ell}^{\ell+n-1}\) blue, a contradiction. Let \(x\) be the number of local extrema of \(\delta(C)\). Then
\[
(x+1)(n-1)-x\ge q-1\ge kn-1=(2h+1)n-1.
\]
Thus \(\delta(C)\) has at least \(2h\) local extrema. On the other hand, if \(\delta'_j\) is a local maximum, then \(\delta'_{j-1}<\delta'_j>\delta'_{j+1}\), which makes \((c_j,c_{j+1})\) a pair. Since \(C\) contains at most \(h-1\) pairs, there are at most \(h-1\) local maxima. Hence there are at most \(h\) local minima, and consequently at most \(2h-1\) local extrema in total, a contradiction.
\medskip

\noindent\textbf{The local red bound.}
It remains to prove that every $(k+1)$-set contains at most three red edges. We use the notation $X+x=X\cup\{x\}$. We first separate the two types.

\begin{claim}\label{no two type}
No $(k+1)$-set contains both a structural $k$-subset and a gap-monotone $k$-subset.
\end{claim}

\begin{proof}[Proof of Claim \ref{no two type}]
Suppose, for contradiction, that $e=(x_1,\ldots,x_k)$ is structural, with $x_1<\cdots<x_k$, and that $e'=e-x_i+x'$ is gap-monotone. Observe that $e$ has $h$ pairs: under $(\Xi_2)$ these are $(x_2,x_3),\allowbreak\ldots,\allowbreak(x_{2h},x_{2h+1})$, while under $(\Xi_3)$ they are $(x_1,x_2),\allowbreak\ldots,\allowbreak(x_{2h-1},x_{2h})$. Put $A=e\cap e'=e-x_i$. Deleting $x_i$ destroys at most one of these pairs, so Fact~\ref{three properties}(1) gives at least $h-1\ge2$ pairs in $A$. On the other hand, $A\subseteq e'$, and Property~III implies that $\delta(A)$ is monotone. Such a set has exactly one pair, a contradiction.
\end{proof}
Suppose, for a contradiction, that $Y=\{y_1<\cdots<y_{k+1}\}$ contains at least four red edges. By Claim \ref{no two type}, they are either all gap-monotone or all structural. Put $\delta_i=\delta(y_i,y_{i+1})$ for $i\in[k]$.

\medskip
\noindent\emph{Case 1: All four red edges are gap-monotone.}
They are either all gap-increasing or all gap-decreasing. Indeed, a gap-increasing and a gap-decreasing edge have at least three vertices in common since $k\ge 7$, which is impossible.

Suppose first that they are all gap-increasing. If $\delta_j>\delta_{j+1}$ for some $j$, every gap-increasing $k$-subset must omit one of $y_j,y_{j+1},y_{j+2}$, since otherwise its restriction to that triple would be decreasing. There are at most three such $k$-subsets, a contradiction. Thus $\delta_1<\delta_2<\cdots<\delta_k$. By Property II, for $i\le k$, we have $\delta(Y-y_i)=\delta(Y)-\delta_i$ and $\delta(Y-y_{k})=\delta(Y-y_{k+1})$. Hence these $4$ red edges give rise to at least $3$ red edges in $K^{(k-1)}_N$, a contradiction. In the gap-decreasing case the same argument applies.
\medskip

\noindent\emph{Case 2: All four red edges are structural.}
Let $e_i=Y-y_i$ for $i\in[k+1]$ and $$\mathcal S=\{i\in[k+1]:e_i~\text{is}~\text{structural}\}.$$
Let $\delta_p$ be the unique minimum of $\delta(Y)$ for some $p\in[k]$. The uniqueness of $\delta_{p}$ follows from Property II.
Let
\[
 L=\{y_1,\ldots,y_p\},\qquad
 R=\{y_{p+1},\ldots,y_{k+1}\}.
\]
For $2\le p\le k-1$, both $L\cap e_i$ and $R\cap e_i$ are nonempty.
Writing $\delta(e_i)=\{\delta^{(i)}_1,\ldots,\delta^{(i)}_{k-1}\}$, Property II gives
\[
 \min_{t\in[k-1]}\delta^{(i)}_t=\delta_p
 =\begin{cases}
\delta_{p-1}^{(i)}, & i\le p,\\
\delta_p^{(i)}, & i>p.
\end{cases}
\]
Consequently, by $(\Xi_2)$ and $(\Xi_3)$,
\begin{equation}\label{case2:min-position}
 i\in\mathcal S
 \quad\Longrightarrow\quad
 \begin{cases}
p-1\in\{1,k-1\}, &i\le p,\\
p\in\{1,k-1\}, &i>p.
\end{cases}
\end{equation}
If $3\le p\le k-2$, then ~\eqref{case2:min-position} yields $\mathcal S=\varnothing$, a contradiction. If $p\in\{2,k-1\}$, then again by~\eqref{case2:min-position},
\[
 \begin{aligned}
 p=2&\quad\Longrightarrow\quad\mathcal S\subseteq\{1,2\},\\
 p=k-1&\quad\Longrightarrow\quad\mathcal S\subseteq\{k,k+1\}.
 \end{aligned}
\]
Hence $|\mathcal S|\le2<4$, a contradiction.

Now, we consider $p=1$. Let $Q=Y-y_1$ and $\delta(Q-y_j)=\{b^{(j)}_1,\ldots,b^{(j)}_{k-2}\}$ for $j\in[2,k+1]$.
Property II gives 
\begin{equation}\label{case2:deletion}
 (b^{(j)}_1,\ldots,b^{(j)}_{k-2})=
 \begin{cases}
 (\delta_3,\ldots,\delta_{k}),&j=2,\\
 (\delta_2,\ldots,\delta_{j-2},\min\{\delta_{j-1},\delta_j\},\delta_{j+1},\ldots,\delta_{k}),
       &3\le j\le k,\\
 (\delta_2,\ldots,\delta_{k-1}),&j=k+1.
 \end{cases}
\end{equation}
Let
\[
 I=\{j\in[2,k+1]:b^{(j)}_1>b^{(j)}_2<b^{(j)}_3>\cdots<b^{(j)}_{k-2}\}.
\]
Since $\delta(y_1,y_j)=\min_{i\in[j-1]}\delta_i=\delta_p=\delta_1$ for all $j\in[2,k+1]$, $(\Xi_2)$ yields
\begin{equation}\label{case2:count}
 Y-y_j\text{ is structural}\ \Longleftrightarrow\ j\in I.
\end{equation}
Thus, 
\[|\mathcal S|=
\begin{cases}
|I|+1, & Q\text{ is structural},\\
|I|, & Q\text{ is not structural}.
\end{cases}\]

\begin{claim}\label{ile 3}
    $|I|\le3$. In particular, if $Q$ is structural, then $|I|=1$.
\end{claim}

\begin{proof}[Proof of Claim \ref{ile 3}]
Suppose, to the contrary, that $|I|\ge4$. Then there exist $j_1,j_2\in I$ with $j_1<j_2$ and $j_2-j_1\ge3$.

If $j_1$ is even, then $2\le j_1\le k-3$. By~\eqref{case2:deletion} and the zigzag inequalities,
\begin{equation*}
\delta_{j_1+1}
=b^{(j_2)}_{j_1}
<b^{(j_2)}_{j_1+1}
\le\delta_{j_1+2}.
\end{equation*}
Here $b^{(j_2)}_{j_1+1}=\min\{\delta_{j_1+2},\delta_{j_1+3}\}$ when $j_2=j_1+3$, and $b^{(j_2)}_{j_1+1}=\delta_{j_1+2}$ otherwise. On the other hand,
\begin{equation*}
\delta_{j_1+1}
=b^{(j_1)}_{j_1-1}
>b^{(j_1)}_{j_1}
=\delta_{j_1+2},
\end{equation*}
a contradiction.

If $j_1$ is odd, then $3\le j_1\le k-2$. Again,~\eqref{case2:deletion} gives
\begin{equation*}
\delta_{j_1}
=b^{(j_2)}_{j_1-1}
<b^{(j_2)}_{j_1}
=\delta_{j_1+1},
\end{equation*}
whereas
\begin{equation*}
\min\{\delta_{j_1-1},\delta_{j_1}\}
=b^{(j_1)}_{j_1-2}
>b^{(j_1)}_{j_1-1}
=\delta_{j_1+1}.
\end{equation*}
This contradicts $\min\{\delta_{j_1-1},\delta_{j_1}\}\le\delta_{j_1}$. Therefore $|I|\le3$.

Now suppose that $Q$ is structural. If $Q$ satisfies $(\Xi_2)$, then $\delta_2=\min_{i\in[2,k]}\delta_i$ is its unique minimum and $\delta_3,\ldots,\delta_k$ is zigzag, so $2\in I$. For every $j\in[3,k+1]$, however,~\eqref{case2:deletion} gives $b^{(j)}_1=\delta_2<b^{(j)}_2$, so $j\notin I$. Thus $I=\{2\}$.

If $Q$ satisfies $(\Xi_3)$, then $\delta_k=\min_{i\in[2,k]}\delta_i$ is its unique minimum and $\delta_2,\ldots,\delta_{k-1}$ is zigzag, so $k+1\in I$. For every $j\in[2,k]$, we have $b^{(j)}_{k-3}>b^{(j)}_{k-2}=\delta_k$, contrary to the final zigzag inequality $b^{(j)}_{k-3}<b^{(j)}_{k-2}$. Hence $I=\{k+1\}$. In either case, $|I|=1$.
\end{proof}

It follows from Claim \ref{ile 3} that
\[
 |\mathcal S|\le
 \begin{cases}
 3,&Q\text{ is not structural},\\
 2,&Q\text{ is structural}.
 \end{cases}
\]
This contradicts $|\mathcal S|\ge4$.

Finally, suppose that $p=k$. For a binary vertex $x$, let $\overline{x}$ denote its coordinatewise complement, and put $z_i=\overline{y_{k+2-i}}$ for $i\in[k+1]$, and set $Z=\{z_1,\ldots,z_{k+1}\}$. Then $z_1<\cdots<z_{k+1}$ and
\begin{equation*}
\delta(z_i,z_{i+1})=\delta_{k+1-i}
\qquad(i\in[k]).
\end{equation*}
Thus $\delta(z_1,z_2)=\min_{i\in[k]}\delta(z_i,z_{i+1})$. Since $k-2$ is odd, reversing a zigzag sequence of length $k-2$ preserves its zigzag inequalities. Consequently, coordinatewise complementation interchanges $(\Xi_2)$ and $(\Xi_3)$ and preserves the number of structural $k$-subsets. Applying the case $p=1$ to $Z$ therefore gives $|\mathcal S|\le3$, a contradiction. This completes the proof of Case 2.
\medskip

We have proved that $\chi$ has at most three red edges on every $(k+1)$-set and no blue $2kn$-set. Finally, take
\[
 M=r_{k-1}(k,3;n)-1
\]
and choose $\varphi$ witnessing this Ramsey lower bound. The construction gives $2^M$ vertices, and hence
\[
 r_k(k+1,4;2kn)>2^M
 =2^{\,r_{k-1}(k,3;n)-1}.
\]
This completes the proof of Lemma \ref{su:oddstep}.
\end{proof}

At the end of this subsection, we prove Theorem~\ref{su:tower}. By (\ref{r_k(k+1,t;n)}) and (\ref{mubayi-suk results}), it suffices to show that
\[
r_k(k+1,t;n)>\operatorname{twr}_{t-1}(c_{k,t}n^{k-t+1}\log n)
\]
for \(4\le t\le k-2\) with \(k-t\) odd, and, to achieve this, we will use Lemma~\ref{su:oddstep} and the following theorem, which is due to Mubayi and Suk \cite[Theorem 7]{MSEH}.
\begin{theorem}[Mubayi and Suk \cite{MSEH}]\label{Mubayianssuk}
    Let $k\ge 6$ and $t\ge 4$. If we are not in the case when $t=4$ and $k$ is odd, then $r_k(k+1,t;2kn)>2^{r_{k-1}(k,t-1;n)-1}$.
\end{theorem}

\begin{proof}[Proof of Theorem~\ref{su:tower} when $k-t$ is odd]
Put $s=k-t+3$ and $B=\prod_{a=s+1}^{k}(2a)$. For the duration of this proof, write $R_{a,b}(u)=r_a(a+1,b;u)-1$. Since Ramsey numbers are integers, Lemma~\ref{su:oddstep} and Theorem~\ref{Mubayianssuk} give, in their respective ranges,
\begin{equation*}
R_{a,b}(2au)\ge 2^{R_{a-1,b-1}(u)}.
\end{equation*}
The first lift has $a=s+1=k-t+4\ge7$ odd and $b=4$, so Lemma~\ref{su:oddstep} applies. The remaining $t-4$ lifts have $b\ge5$ and are covered by Theorem~\ref{Mubayianssuk}. Thus, for all sufficiently large integers $m$,
\begin{equation}\label{last equ}
R_{k,t}(Bm)\ge\operatorname{twr}_{t-2}\bigl(R_{s,3}(m)\bigr).
\end{equation}
After decreasing a positive constant $c=c_{k,t}$, Theorem~\ref{thm:localbase} gives $R_{s,3}(m)\ge2^{c m^{s-2}\log m}$ for all sufficiently large $m$. Since $s-2=k-t+1$, equation~\eqref{last equ} yields
\begin{equation*}
R_{k,t}(Bm)\ge\operatorname{twr}_{t-1}\bigl(c m^{k-t+1}\log m\bigr).
\end{equation*}
For arbitrary sufficiently large $n$, take $m=\lfloor n/B\rfloor$. Then $Bm\le n$ and $m=\Theta_{k,t}(n)$. Monotonicity and a further decrease of the constant give
\begin{equation*}
r_k(k+1,t;n)>\operatorname{twr}_{t-1}\bigl(c_{k,t}n^{k-t+1}\log n\bigr).
\end{equation*}\end{proof}
This completes the proof of Theorem \ref{su:tower}.
 
\section*{Remarks}

\noindent\textbf{Declaration on the use of AI:}~~The authors acknowledge the use of AI tools during the preparation of this work. The research initially focused on the Ramsey numbers \(r_k(B^{(k)}, n)\). With input from ChatGPT-6, the authors obtained an initial proof of the desired lower bound. However, this initial proof was technically complicated. The authors subsequently developed a new approach, leading to a simpler and more accessible proof.

After writing the manuscript, the authors recognized the potential for further generalizing the method and engaged in several rounds of discussions with ChatGPT-6 on this question. Then the authors extended the method to obtain lower bounds for \(r_k(\mathcal{F}_{k,r} \cup \mathcal{H}_{k,r}, n)\), yielding further results on \(r_k(k+1,3; n)\) and \(r_k(F_{\ell,k,d}, n)\). Eventually, the authors independently developed a new stepping-up method, Lemma~\ref{su:oddstep}, which extends the results on \(r_k(k+1,3; n)\) to \(r_k(k+1,t; n)\).

The final manuscript was written by the authors. The authors take full responsibility for all AI-assisted aspects in this paper.


\medskip

\noindent\textbf{Note added:}~~After completing the manuscript, the authors learned from Ruilong Liu, who was visiting the Georgia Institute of Technology, that Eion Mulrenin and Cosmin Pohoata had independently obtained partially overlapping results.

\end{spacing}
\end{document}